\documentclass[11pt]{article}

\usepackage[margin=1in]{geometry}
\usepackage{bm,amsmath,commath,mathtools,amsfonts,amssymb}
\usepackage{graphicx, caption, subcaption, float}
\usepackage{url,cite,hyperref}
\hypersetup{
    colorlinks=true,
    linkcolor=blue,
    citecolor = red,
}

\usepackage{verbatim}

\usepackage{authblk}
\usepackage{onimage}

\newcommand{\vect}[1]{\bm{#1}}
\newcommand{\vF}{\vect{F}}

\newcommand{\vu}{\vect{u}}
\newcommand{\vv}{\vect{v}}

\newcommand{\vx}{\vect{x}}
\newcommand{\vy}{\vect{y}}

\newcommand{\mat}[1]{\bm{#1}}
\newcommand{\mA}{\mat{A}}
\newcommand{\mB}{\mat{B}}
\newcommand{\mC}{\mat{C}}
\newcommand{\mI}{\mat{I}}
\newcommand{\mP}{\mat{P}}
\newcommand{\mQ}{\mat{Q}}

\newcommand{\mU}{\mat{U}}
\newcommand{\mV}{\mat{V}}
\newcommand{\mX}{\mat{X}}
\newcommand{\mY}{\mat{Y}}

\newcommand{\mGamma}{\mat{\Gamma}}

\DeclareMathOperator{\rank}{Rank}
\DeclareMathOperator{\imag}{Im}

\title{Online dynamic mode decomposition for time-varying systems}

\author[1]{Hao Zhang\footnote{Corresponding author. Email: haozhang@princeton.edu}}
\author[1]{Clarence W. Rowley}
\author[2]{Eric A. Deem}
\author[2]{Louis N. Cattafesta}
\affil[1]{Mechanical and Aerospace Engineering, Princeton University}
\affil[2]{Mechanical Engineering, Florida State University}

\date{\today}

\begin{document}
\maketitle

\begin{abstract}
  Dynamic mode decomposition (DMD) is a popular technique for modal
  decomposition, flow analysis, and reduced-order modeling.  In situations where
  a system is time varying, one would like to update the
  system's description online as time evolves. This work provides an efficient
  method for computing DMD in real time, updating the approximation of a
  system's dynamics as new data becomes available.  The algorithm does not
  require storage of past data, and computes the exact DMD matrix using rank-1
  updates.  A weighting factor that places less weight on older data can be
  incorporated in a straightforward manner, making the method particularly well
  suited to time-varying systems.  A variant of the method may also be applied
  to online computation of ``windowed DMD'', in which only the most recent data
  are used.  The efficiency of the method is compared against several existing
  DMD algorithms: for problems in which the state dimension is less than
  about~200,
  the proposed algorithm is the most efficient for real-time
  computation, and it can be orders of magnitude more efficient than the
  standard DMD algorithm. The method is demonstrated on several examples,
  including a time-varying linear system and a more
  complex example using data from a wind tunnel experiment.  In particular, we
  show that the method is effective at capturing the dynamics of surface
  pressure measurements in the flow over a flat plate with an unsteady
  separation bubble.
\end{abstract}

\section{Introduction}
Modal decomposition methods are widely used in studying complex dynamical systems
such as fluid flows. In particular, dynamic mode decomposition (DMD)
\cite{schmid2010dynamic,rowley2009spectral} has become increasingly popular in
the fluids community. DMD decomposes spatio-temporal data into spatial modes
(DMD modes) each of which has simple temporal behavior characterized by single
frequency and growth/decay rate (DMD eigenvalues). DMD has been successfully
applied to a wide range of problems, for instance as discussed
in~\cite{rowley2017model, kutz2016dynamic}. The idea of DMD is to fit a linear
system to observed dynamics.  However, DMD is also a promising technique for
nonlinear systems, as it has been shown to be a finite-dimensional approximation
to the Koopman operator, an infinite-dimensional linear operator that captures
the full behavior of a nonlinear dynamical system
\cite{rowley2009spectral,tu2014dynamic}.

Recently, several algorithms have been proposed to compute DMD modes efficiently
for very large datasets, for instance using randomized methods
\cite{erichson2016randomized,erichson2017randomized}. In situations in which
the incoming data is ``streaming'' in nature, and one does not wish to store all
of the data, a ``streaming DMD'' algorithm performs online updating of the DMD
modes and eigenvalues \cite{hemati2014dynamic}.  Streaming DMD keeps track of a
small number of orthogonal basis vectors and
updates the DMD matrix projected onto the corresponding subspace.
Another related method uses an incremental SVD algorithm to compute DMD modes on the fly \cite{Matsumoto:2017}.
The work
proposed here may be viewed as an alternative to Streaming DMD, in that we
provide a method for updating the DMD matrix in real time, without the need to
store all the raw data.  Our method
differs from Streaming DMD in that we compute the exact DMD matrix, rather than
a projection onto basis functions; in addition, we propose various methods for
better approximation of time-varying dynamics, in particular by ``forgetting''
older snapshots, or giving them less weight than more recent snapshots.

The paper is organized as follows.  In Section~\ref{sec:onlineDMD}, we give an
overview of DMD, and describe the online DMD algorithm.  In
Section~\ref{sec:windowDMD}, we discussed a variant called windowed DMD, in
which only the most recent data are used.   In Section~\ref{sec:system-id}, we
briefly describe how these methods may be used in online system identification,
and in Section~\ref{sec:result} we compare the different algorithms on various
examples.

\section{Online dynamic mode decomposition}
\label{sec:onlineDMD}
\subsection{The problem}
\label{sec:dmd-problem}

We first give a brief summary of the standard DMD algorithm, as described
in~\cite{tu2014dynamic}.  Suppose we have a discrete-time dynamical system given
by
\begin{equation*}
\vx_{j+1} = \vF(\vx_j),
\end{equation*}
where $\vx_j \in \mathbb{R}^n$ is the state vector, and
$\vF:\mathbb{R}^n\to\mathbb{R}^n$ defines the dynamics.
For a given state $\vx_j$, let $\vy_j=\vF(\vx_j)$; we call $(\vx_j,\vy_j)$ a
{\em snapshot pair}.  For DMD, we assume we have access to a collection of
snapshot pairs $(\vx_j,\vy_j)$, for $j=1,\ldots,k$.  (It is often the case that
$\vx_{j+1}=\vy_j$, corresponding to a sequence of points along a single
trajectory, but this is not required.)

DMD seeks to find a matrix $\mA$ such that $\vy_j=\mA\vx_j$, in an approximate
sense.  DMD modes are then eigenvectors of the matrix~$\mA$, and DMD eigenvalues
are the corresponding eigenvalues.  In the present work, we are interested in
obtaining a matrix $\mA$ that varies in time, giving us a local linear model for
the dynamics, but in the standard DMD approach, one seeks a single matrix~$\mA$.

Given snapshot pairs $(\vx_j,\vy_j)$ for $j=1,\ldots,k$, we form matrices
\begin{equation}
\label{eq:5}
\mX_k = \begin{bmatrix} \vx_{1} & \vx_{2} & \cdots & \vx_{k}\end{bmatrix},\qquad
\mY_k = \begin{bmatrix} \vy_{1} & \vy_{2} & \cdots & \vy_{k}\end{bmatrix},
\end{equation}
which both have dimension $n\times k$.  We wish to find an $n\times n$ matrix
$\mA_k$ such that $\mA_k \mX_k = \mY_k$ approximately holds; in particular, we
are interested in the {\em overconstrained\/} problem, in which $k>n$. When the
problem is {\em underconstrained\/}, the model will tend to overfit the data,
and any noise present in the data will lead to poor performance of the model
\cite{hawkins2004problem}. The DMD
matrix $\mA_k$ is then found by minimizing the cost function
\cite{schmid2010dynamic,rowley2009spectral}
\begin{equation} \label{eqn:onlinecost}
J_k = \sum_{i=1}^k \|\vy_i - \mA_k \vx_i\|^2 = \|\mY_k - \mA_k \mX_k\|_F^2,
\end{equation}
where $\|\cdot\|$ denotes the Euclidean norm on vectors and $\|\cdot\|_F$
denotes the Frobenius norm on matrices.  The unique minimum-norm solution to this least-squares problem is given by
\begin{equation}\label{eqn:onlineA}
\mA_k = \mY_k \mX^{+}_k,
\end{equation}
where $\mX^{+}_k$ denotes the Moore-Penrose pseudoinverse of $\mX_k$.

Here, we shall assume that $\mX_k$ has full row rank, in which case
$\mX_k\mX_k^T$ is invertible, and
\begin{equation}
  \label{eq:2}
  \mX_k^+=\mX_k^T(\mX_k\mX_k^T)^{-1}.
\end{equation}
This assumption is essential for the development of the online algorithm, as we
shall see shortly. Under this assumption, the $\mA_k$ given above is the unique
solution that minimizes $J_k$.  In this paper, we are interested in the case in
which the number of snapshots $k$ is large, compared with the state
dimension~$n$.

Our primary focus here is systems that may be slowly varying in time, so that
the matrix $\mA_k$ should evolve as $k$ increases.  In the following section, we
will present an efficient algorithm for updating $\mA_k$ as more data becomes
available.  Furthermore, if the system is time varying, it may make sense to
weight more recent snapshots more heavily than less recent snapshots.  In this
spirit, we will consider minimizing a modified cost function
\begin{equation}
  \label{eq:1}
  \tilde J_k = \sum_{i=1}^k \rho^{k-i} \|\vy_i - \mA_k \vx_i\|^2,
\end{equation}
for some constant $\rho$ with $0<\rho\le 1$.  When $\rho=1$, this cost function
is the same as~\eqref{eqn:onlinecost}, and when $\rho<1$, errors in past
snapshots are discounted.   Our algorithm will apply to this minimization
problem as well, with only minor modifications and no increase in computational
effort.


\begin{figure}[!htb]
        \centering
     \includegraphics[width=.7\linewidth]{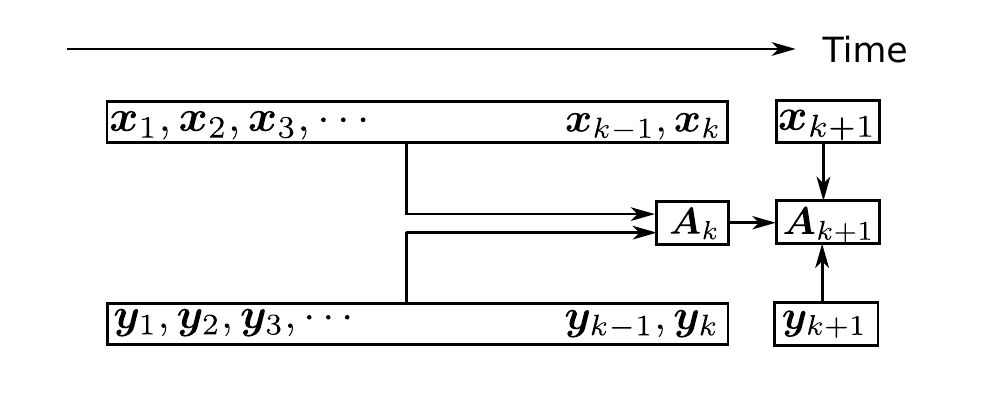}
\caption{A cartoon of the online DMD setup. $\mA_k$ is the optimal (least-squares) fit that maps $\mX_k=[\vx_{1},\vx_{2},\cdots,\vx_{k}]$ to $\mY_k=[\vy_{1},\vy_{2},\cdots,\vy_{k}]$. Arrow indicates the information flow, and box denotes block of information. At time $k+1$, $\mA_{k}$ is updated to find $\mA_{k+1}$, using the information from time $k$, and new available snapshot pair $\vx_{k+1},\vy_{k+1}$ at time $k+1$. $\mA_{k+1}$ is the optimal (least-squares) fit that maps $\mX_{k+1}=[\vx_{1},\vx_{2},\cdots,\vx_{k},\vx_{k+1}]$ to $\mY_{k+1}=[\vy_{1},\vy_{2},\cdots,\vy_{k},\vy_{k+1}]$.}
\label{fig:onlineDMD}
\end{figure}

A sketch of the online DMD setup is shown in Figure \ref{fig:onlineDMD}.
Suppose we have already computed $\mA_k$ for a given dataset.  As time
progresses and a new pair of snapshots $(\vx_{k+1},\vy_{k+1})$ becomes
available, the matrix $\mA_{k+1}$ may be updated according to the formula given
in equation (\ref{eqn:onlineA}). If $\mA_{k+1}$ is computed directly in this
manner, we call this the ``standard approach''.

There are two drawbacks to the ``standard approach''.   First, it requires
computing the pseudoinverse of $\mX_k$ whenever new snapshots are acquired, and
for this reason it is computationally expensive.  In addition, the method
requires storing all the snapshots (i.e., storing the matrix $\mX_k$), which may
be challenging or impossible as the number of snapshots $k$ increases.

 \subsection{Algorithm for online DMD}
\label{sec:alg-online-dmd}
To overcome the above two shortcomings, we propose a different approach to find
the solution to equation (\ref{eqn:onlineA}) in the ``online setting'', in which
we want to compute $\mA_{k+1}$ given a matrix $\mA_{k}$ and a new pair of
snapshots $(\vx_{k+1},\vy_{k+1})$.  The algorithm we present is based on the
idea that $\mA_{k+1}$ should be close to $\mA_k$ in some sense.


First, observe that, using~\eqref{eq:2}, we may write~\eqref{eqn:onlineA} as
\begin{equation} \label{eqn:onlineformAk}
\mA_k = \mY_k \mX_k^T (\mX_k \mX_k^T)^{-1}= \mQ_k \mP_k,
\end{equation}
where $\mQ_k$ and $\mP_k$ are $n\times n$ matrices given by
\begin{subequations}
\label{eq:3}
\begin{align}
\mQ_k &= \mY_k \mX_k^T,\\
\mP_k &= (\mX_k \mX_k^T)^{-1}.\label{eqn:onlinedefPk}
\end{align}
\end{subequations}
The condition that $\mX_k$ has rank~$n$ ensures that $\mX_k \mX_k^T$ is
invertible, and hence $\mP_k$ is well defined. Note also that $\mP_k$ is
symmetric and strictly positive definite.

At time $k+1$, we wish to compute $\mA_{k+1} = \mQ_{k+1} \mP_{k+1}$. Clearly, $\mQ_{k+1}, \mP_{k+1}$ are related to $\mQ_{k}, \mP_{k}$:
\begin{align*}
\mQ_{k+1} = \mY_{k+1} \mX_{k+1}^T =
\begin{bmatrix}
\mY_{k} & \vy_{k+1}
\end{bmatrix}
\begin{bmatrix}
\mX_{k} & \vx_{k+1}
\end{bmatrix}^T
&= \mY_k \mX_k^T + \vy_{k+1} \vx_{k+1}^T,\\
\mP_{k+1}^{-1}  = \mX_{k+1} \mX_{k+1}^T =
\begin{bmatrix}
\mX_{k} & \vx_{k+1}
\end{bmatrix}
\begin{bmatrix}
\mX_{k} & \vx_{k+1}
\end{bmatrix}^T
&= \mX_k \mX_k^T + \vx_{k+1} \vx_{k+1}^T.
\end{align*}
Because $\mX_k$ already has rank~$n$, and adding an additional column
cannot reduce the rank of a matrix, it follows that $\mX_{k+1}$ also has
rank~$n$, so $\mP_{k+1}$ is well defined.  The above equation shows that, given
$\mQ_k$ and~$\mP_k^{-1}$, we may find $\mQ_{k+1}$ and~$\mP_{k+1}^{-1}$ with
simple rank-1 updates:
\begin{align*}
\mQ_{k+1} &= \mQ_k + \vy_{k+1} \vx_{k+1}^T,\\
\mP_{k+1}^{-1} &=  \mP_k^{-1} + \vx_{k+1} \vx_{k+1}^T. \label{eqn:onlinedefPk1}
\end{align*}
The updated DMD matrix is then given by
\begin{equation} \label{eqn:onlinedefAk1}
\mA_{k+1} = \mQ_{k+1} \mP_{k+1} = (\mQ_k + \vy_{k+1} \vx_{k+1}^T) (\mP_k^{-1} + \vx_{k+1} \vx_{k+1}^T)^{-1}.
\end{equation}

Then the problem is reduced to how to find $\mP_{k+1}$ from
$\mP_{k}$ in an efficient manner. Computing the inverse directly would require
$\mathcal{O}(n^3)$ operations, and would not be efficient.  However, because
$\mP_{k+1}$ is the inverse of a rank-1 update of $\mP_{k}^{-1}$, we may take
advantage of a matrix inversion formula known as the Sherman-Morrison formula \cite{sherman1950adjustment, hager1989updating}.

Suppose $\mA$ is an invertible square matrix, and $\vu, \vv$ are column
vectors. Then $\mA + \vu\vv^T$ is invertible if and only if $1+\vv^T\mA\vu \neq
0$, and in this case, the inverse is given by the Sherman-Morrison formula
\begin{equation}
\label{eq:Sherman-Morrison}
(\mA + \vu\vv^T)^{-1} = \mA^{-1} - \frac{\mA^{-1}\vu\vv^T\mA^{-1}}{1+\vv^T \mA^{-1} \vu}.
\end{equation}
This formula is a special case of the more general matrix
inversion lemma (or Woodbury formula) \cite{woodbury1950inverting, hager1989updating}.

Applying the formula to the expression for $\mP_{k+1}$, we obtain
\begin{subequations}
\label{eqn:onlineupdatePk1}
\begin{equation}
\mP_{k+1} = (\mP_k^{-1} + \vx_{k+1} \vx_{k+1}^T)^{-1} = \mP_k - \gamma_{k+1} \mP_k \vx_{k+1}\vx_{k+1}^T\mP_k,
\end{equation}
where
\begin{equation} \label{eqn:gammak1}
\gamma_{k+1} = \frac{1}{1+\vx_{k+1}^T\mP_k\vx_{k+1}}.
\end{equation}
\end{subequations}
Note that, because $\mP_k$ is positive definite, the scalar quantity
$1+\vx_{k+1}^T\mP_k\vx_{k+1}$ is always nonzero, so the formula applies.
Therefore, the updated DMD matrix may be written
\begin{equation} \label{eqn:onlineexpandAk1}
\begin{split}
\mA_{k+1} & = (\mQ_k + \vy_{k+1} \vx_{k+1}^T) (\mP_k - \gamma_{k+1} \mP_k \vx_{k+1}\vx_{k+1}^T\mP_k) \\
& = \mQ_k \mP_k - \gamma_{k+1} \mQ_k \mP_k \vx_{k+1}\vx_{k+1}^T\mP_k \\
& + \vy_{k+1} \vx_{k+1}^T \mP_k - \gamma_{k+1} \vy_{k+1} \vx_{k+1}^T \mP_k \vx_{k+1}\vx_{k+1}^T\mP_k.
\end{split}
\end{equation}
We can simplify the last two terms, since
\begin{equation*}
\begin{split}
\vy_{k+1} \vx_{k+1}^T \mP_k - \gamma_{k+1} \vy_{k+1} \vx_{k+1}^T \mP_k \vx_{k+1}\vx_{k+1}^T\mP_k & = \gamma_{k+1} \vy_{k+1}(\gamma_{k+1}^{-1}-\vx_{k+1}^T \mP_k \vx_{k+1})\vx_{k+1}^T\mP_k \\
& = \gamma_{k+1} \vy_{k+1} \vx_{k+1}^T\mP_k,
\end{split}
\end{equation*}
where we have used~\eqref{eqn:gammak1}.  Substituting into~\eqref{eqn:onlineexpandAk1}, we obtain
\begin{align*}
\mA_{k+1} & = \mQ_k \mP_k - \gamma_{k+1} \mQ_k \mP_k \vx_{k+1}\vx_{k+1}^T\mP_k + \gamma_{k+1} \vy_{k+1} \vx_{k+1}^T\mP_k \\
& = \mA_k - \gamma_{k+1} \mA_k \vx_{k+1}\vx_{k+1}^T\mP_k + \gamma_{k+1} \vy_{k+1} \vx_{k+1}^T\mP_k,
\end{align*}
and hence
\begin{equation}
  \label{eqn:onlineupdateAk1}
  \mA_{k+1} = \mA_k + \gamma_{k+1}(\vy_{k+1}-\mA_k\vx_{k+1})\vx_{k+1}^T\mP_k.
\end{equation}
The above formula gives a rule for computing $\mA_{k+1}$ given $\mA_k, \mP_k$
and the new snapshot pair $(\vx_{k+1}, \vy_{k+1})$. In order to use this formula
recursively, we also need to compute $\mP_{k+1}$ using~\eqref{eqn:onlineupdatePk1}, given $\mP_k$ and $\vx_{k+1}$.

There is an intuitive interpretation for the update
formula~\eqref{eqn:onlineupdateAk1}. The quantity $(\vy_{k+1}-\mA_k\vx_{k+1})$
can be considered as the prediction error from the current model~$\mA_k$, and
the DMD matrix is updated by adding a term proportional to this error.

The updates in equations \eqref{eqn:onlineupdatePk1}
and~\eqref{eqn:onlineupdateAk1} together require only two matrix vector
multiplications ($\mA_k\vx_{k+1}$ and $\mP_k\vx_{k+1}$, since $\mP_k$ is
symmetric), and two vector outer products, for a total of $4n^2$ floating-point
multiplies.
This is
much more efficient than applying the standard DMD algorithm, which involves a
singular value decomposition or pseudoinverse, and requires $O(kn^2)$
multiplies, where $k>n$.  In our approach, two $n \times n$ matrices need to be
stored ($\mA_k$ and~$\mP_k$), but the large $n \times k$ snapshot matrices $(\mX_k,\mY_k)$ do not need to be
stored.

It is worth emphasizing that the update formulas \eqref{eqn:onlineupdatePk1}
and~\eqref{eqn:onlineupdateAk1} compute the DMD matrix
$\mA_{k+1} = \mY_{k+1} \mX_{k+1}^+$ exactly (up to machine precision).  That is,
with exact arithmetic, our formulas give the same results as the standard DMD
algorithm.  The matrix $\mP_k$ does involve ``squaring up'' the matrix $\mX_k$,
which could in principle lead to difficulties with numerical stability for
ill-conditioned problems \cite{yip1986note,allower1990update}.  However, we have
not encountered problems with numerical stability in the examples we have tried
(see Section~\ref{sec:result}).

\paragraph{Initialization} The algorithm described above needs a starting point.
In particular, to apply the updates~\eqref{eqn:onlineupdatePk1}
and~\eqref{eqn:onlineupdateAk1}, one needs the matrices $\mP_k$ and $\mA_k$ at
timestep~$k$.
The initialization technique is similar to the initialization of recursive least-squares estimation described in \cite{hsia1977system}. Two practical approach are discussed below.
The most straightforward way to initialize the algorithm is to
first collect at least $n$ snapshots (more precisely, enough snapshots so that
$\mX_k$ as defined in~\eqref{eq:5} has rank~$n$), and then compute $\mP_k$ and~$\mA_k$ using the standard
DMD algorithm, from~\eqref{eqn:onlineformAk} and~\eqref{eq:3}:
\begin{equation}
  \label{eq:4}
  \mA_k = \mY_k\mX_k^+,\qquad \mP_k = (\mX_k\mX_k^T)^{-1}.
\end{equation}
If for some reason this is not desirable, then an alternative approach is to
initialize $\mA_0$ to a random matrix (e.g., the zero matrix), and set
$\mP_0=\alpha\mI$, where $\alpha$ is a large positive scalar.  Then in the limit
as $\alpha\to\infty$, the matrices $\mP_k,\mA_k$ computed by the updates
\eqref{eqn:onlineupdatePk1} and~\eqref{eqn:onlineupdateAk1} converge to the true
values given by~\eqref{eq:4}.

\paragraph{Multiple snapshots} In our method, the DMD matrix $\mA_k$ gets
updated at every time step when a new snapshot pair becomes available. In
principle, one could update the DMD matrix less frequently (for instance every 10
time steps). The above derivation can be appropriately modified to handle this
case, using the more general Woodbury formula (see equation~\eqref{eq:9}) \cite{woodbury1950inverting,hager1989updating}. However, if $s$ is the number of new snapshots to be
incorporated, the computational cost of a single rank-$s$ update is roughly the
same as applying the rank-1 formula $s$ times, so there does not appear to be a
benefit to incorporating multiple snapshots at once.

\paragraph{Extensions} As is the case for most DMD algorithms (including Streaming DMD), the online DMD algorithm described above applies
more generally to extended DMD (EDMD) \cite{williams2015data}, simply replacing the
state observations $\vx_k, \vy_k$ by the corresponding vectors of
observables. In addition, the algorithm can be used for real-time online system
identification, including both linear and nonlinear system identification, as we
shall discuss in Section~\ref{sec:system-id}.

\paragraph{Summary} To summarize, the algorithm proceeds as follows:
\begin{enumerate}
\item Collect $k$ snapshot pairs $(\vx_j,\vy_j)$, $j=1,\ldots,k$, where $k>n$ is
  large enough so that $\rank \mX_k=n$ (where $\mX_k$ is given by~\eqref{eq:5}).
\item Compute $\mA_k$ and $\mP_k$ from~\eqref{eq:4}.
\item When a new snapshot pair $(\vx_{k+1},\vy_{k+1})$ becomes available, update
  $\mA_k$ and~$\mP_k$ according to~\eqref{eqn:onlineupdateAk1}
  and~\eqref{eqn:onlineupdatePk1}.
\end{enumerate}
Implementations of this algorithm in both Matlab and Python are publicly
available at~\cite{zhang2017odmd}.

\subsection{Weighted online DMD}
\label{sec:weighted-online}
As mentioned previously, the online DMD algorithm described above is ideally
suited to cases for which the system is varying in time, so that we want to
revise our estimate of the DMD matrix $\mA_k$ in real time.  In such a
situation, we might wish to place more weight on recent snapshots, and gradually
``forget'' the older snapshots, by minimizing a cost function of the
form~\eqref{eq:1} instead of the original cost function~\eqref{eqn:onlinecost}.
This weighting scheme is analogous to that used in real-time least-squares
approximation~\cite{hsia1977system}.
This idea may also be used with streaming DMD, and in fact has been considered
before (the conference presentation~\cite{Hemati:2016} implemented such a
``forgetting factor'' with
streaming DMD, although it did not appear in the associated paper).
It turns out that the online DMD algorithm can be adapted to minimize the cost
function~\eqref{eq:1} with only minor modifications to the algorithm.

We now consider the cost function
\begin{equation*}
\tilde J_k = \sum_{i=1}^k \rho^{k-i} \|\vy_i - \mA_k \vx_i\|^2,\qquad 0 < \rho \leq 1,
\end{equation*}
where $\rho$ is the weighting factor.  For instance, if we wish our snapshots to
have a ``half-life'' of $m$ samples, then we could choose $\rho=2^{-1/m}$.
For convenience, let us take
$\rho=\sigma^2$ where $0<\sigma\le 1$, and write the cost function as
\begin{equation*}
\tilde J_k = \sum_{i=1}^k \|\sigma^{k-i} \vy_i - \mA_k \sigma^{k-i} \vx_i\|^2.
\end{equation*}
If we define matrices based on scaled versions of the snapshots, as
\begin{align*}
  \tilde\mX_k
  &=
    \begin{bmatrix}
      \sigma^{k-1} \vx_{1} & \sigma^{k-2} \vx_{2} & \cdots & \vx_{k}
    \end{bmatrix},\\
  \tilde\mY_k
  &=
    \begin{bmatrix}
      \sigma^{k-1} \vy_{1} & \sigma^{k-2} \vy_{2} & \cdots & \vy_{k}
    \end{bmatrix},
\end{align*}
then the cost function can be written as
\begin{equation*}
\tilde J_k = \|\tilde\mY_k - \mA_k \tilde\mX_k\|_F^2.
\end{equation*}
The unique least-squares solution that minimizes this cost function (assuming $\tilde\mX_k$
has full row rank) is given by
\begin{equation*}
\mA_k = \tilde\mY_k \tilde\mX_k^+ = \tilde\mY_k \tilde\mX_k^T (\tilde\mX_k \tilde\mX_k^T)^{-1} = \tilde\mQ_k \tilde\mP_k,
\end{equation*}
where we define
\begin{align*}
\tilde\mQ_k &= \tilde\mY_k \tilde\mX_k^T,\\
\tilde\mP_k &= (\tilde\mX_k \tilde\mX_k^T)^{-1}.
\end{align*}
At step $k+1$, we wish to compute $\mA_{k+1}  = \tilde\mQ_{k+1} \tilde\mP_{k+1}$. We write down $\tilde\mX_{k+1}, \tilde\mY_{k+1}$ explicitly as
\begin{align*}
  \tilde\mX_{k+1}
  &=
    \begin{bmatrix}
      \sigma^{k} \vx_{1} & \sigma^{k-1} \vx_{2} & \cdots & \sigma \vx_{k} & \vx_{k+1}
    \end{bmatrix}
  =
    \begin{bmatrix}
      \sigma \tilde\mX_{k} & \vx_{k+1}
    \end{bmatrix},\\
  \tilde\mY_{k+1}
  &=
    \begin{bmatrix}
      \sigma^{k} \vy_{1} & \sigma^{k-1} \vy_{2} & \cdots & \sigma \vy_{k} & \vy_{k+1}
    \end{bmatrix}
  =
    \begin{bmatrix}
      \sigma \tilde\mY_{k} & \vy_{k+1}
    \end{bmatrix}.
\end{align*}
Therefore, $\tilde\mQ_{k+1}$ can be written
\begin{align*}
\tilde\mQ_{k+1}
  &= \tilde\mY_{k+1} \tilde\mX_{k+1}^T
    =
    \begin{bmatrix}
      \sigma \tilde\mY_{k} & \vy_{k+1}
    \end{bmatrix}
    \begin{bmatrix}
      \sigma \tilde\mX_{k} & \vx_{k+1}
    \end{bmatrix}^T \\
& = \sigma^2 \tilde\mY_{k} \tilde\mX_{k}^T + \vy_{k+1} \vx_{k+1}^T\\
&= \rho \tilde\mQ_{k} + \vy_{k+1} \vx_{k+1}^T,
\end{align*}
and similarly
\begin{equation} \label{eqn:onlinePk1tilde}
\tilde\mP_{k+1}^{-1} =\rho \tilde\mP_{k}^{-1} + \vx_{k+1} \vx_{k+1}^T.
\end{equation}
The updated DMD matrix is then given by
\begin{equation*}
\mA_{k+1} = \tilde\mQ_{k+1} \tilde\mP_{k+1} = (\rho \tilde\mQ_{k} + \vy_{k+1} \vx_{k+1}^T)(\rho \tilde\mP_{k}^{-1} + \vx_{k+1} \vx_{k+1}^T)^{-1}.
\end{equation*}
As before, we can apply the Sherman-Morrison formula~\eqref{eq:Sherman-Morrison} to~\eqref{eqn:onlinePk1tilde} and obtain
\begin{equation*}
\tilde\mP_{k+1} = \frac{\tilde\mP_k}{\rho} - \gamma_{k+1} \frac{\tilde\mP_k}{\rho} \vx_{k+1} \vx_{k+1}^T \frac{\tilde\mP_k}{\rho},
\end{equation*}
where
\begin{equation*}
\gamma_{k+1} = \frac{1}{1+  \vx_{k+1}^T (\tilde\mP_{k}/\rho) \vx_{k+1}}.
\end{equation*}
Let us rescale $\tilde\mP_k$, and define
\begin{equation*}
\hat\mP_k = \frac{\tilde\mP_k}{\rho} = \frac{1}{\rho} (\tilde\mX_k \tilde\mX_k^T)^{-1}.
\end{equation*}
Then after some manipulation, the formula for $\mA_{k+1}$ becomes
\begin{equation}
\label{eq:6}
\mA_{k+1} = \mA_k + \gamma_{k+1} (\vy_{k+1} - \mA_k \vx_{k+1}) \vx_{k+1}^T \hat\mP_k,
\end{equation}
where
\begin{subequations}
\label{eq:weighed-P-update}
\begin{align}
\hat\mP_{k+1} &= \frac{1}{\rho} (\hat\mP_k - \gamma_{k+1} \hat\mP_k \vx_{k+1} \vx_{k+1}^T \hat\mP_k),\\
\gamma_{k+1} &= \frac{1}{1+  \vx_{k+1}^T \hat\mP_k \vx_{k+1}}.
\end{align}
\end{subequations}
Observe that the update~\eqref{eq:6} for~$\mA_{k+1}$ is identical to the
update~\eqref{eqn:onlineupdateAk1} from the previous section, with $\mP_k$
replaced by $\hat\mP_k$, and the update rule~\eqref{eq:weighed-P-update} for
$\hat\mP_{k+1}$ differs from~\eqref{eqn:onlineupdatePk1} only by a factor
of~$\rho$.  When $\rho=1$, of course, the above formulas are identical to
those given in Section~\ref{sec:alg-online-dmd}.

\section{Windowed dynamic mode decomposition}
\label{sec:windowDMD}

In Section~\ref{sec:weighted-online}, we presented a method for gradually
``forgetting'' older snapshots, by giving them less weight in a cost function.
In this section, we discuss an alternative method, which uses a hard cutoff: in
particular, we consider a ``window'' containing only the most recent snapshots,
for instance as used in~\cite{Grosek:2014}.


\subsection{The problem}
If the dynamics are slowly varying with time, we may wish to use only the most
recent snapshots to identify the dynamics. Here, we consider the case where we
use only a fixed ``window'' containing the most recent snapshots.  Here, we
present an ``online'' algorithm to compute windowed DMD efficiently, again using
low-rank updates, as in the previous section.  We refer to the resulting method
as ``windowed dynamic mode decomposition'' (windowed DMD).

At time $t_k$, suppose we have access to past snapshot pairs $\{(\vx_j,\vy_j)\}_{j=k-w+1}^k$ in a finite time
window of size $w$.  We would like to fit a linear
model $\mA_k$, such that $\vy_j=\mA_k\vx_j$ (at least approximately) for all $j$
in this window.  Let
\begin{equation}
\label{eqn:windowXk}
\mX_k = \begin{bmatrix} \vx_{k-w+1} & \vx_{k-w+2} & \cdots & \vx_{k}\end{bmatrix},\qquad
\mY_k = \begin{bmatrix} \vy_{k-w+1} & \vy_{k-w+2} & \cdots & \vy_{k}\end{bmatrix},
\end{equation}
both $n\times w$ matrices.  Then we seek an $n\times n$  matrix $\mA_k$ such
that $\mA_k \mX_k = \mY_k$ approximately holds. More precisely (as explained in Section~\ref{sec:dmd-problem}), the DMD matrix
$\mA_k$ is found by minimizing
\begin{equation} \label{eqn:windowcost}
J_k = \|\mY_k - \mA_k \mX_k\|_F^2.
\end{equation}
As before, we assume that the rank of $\mX_k$ is~$n\le w$, so that there is a unique
solution to this least-squares problem, given by
\begin{equation}\label{eqn:windowA}
\mA_k = \mY_k \mX^{+}_k,
\end{equation}
where $\mX^{+}_k=\mX_k^T(\mX_k\mX_k^T)^{-1}$ is the Moore-Penrose pseudoinverse
of $\mX_k$.  (Note, in particular, that we require that the window size~$w$ be
at least as large as the state dimension~$n$, so that $\mX_k\mX_k^T$ is invertible.)

\begin{figure}[!htb]
        \centering
     \includegraphics[width=.7\linewidth]{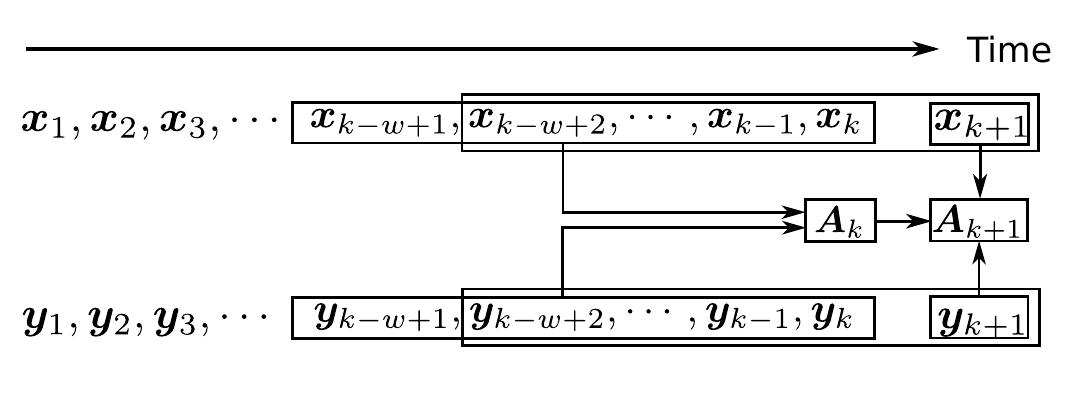}
\caption{A cartoon of the windowed DMD setup. At time $k$, $\mA_k$ depends only
  on the $w$ most recent snapshots.  At time $k+1$, one new snapshot is added, and
  the oldest snapshot is dropped.}
\label{fig:windowDMD}
\end{figure}

A sketch of the windowed DMD setup is shown in Figure \ref{fig:windowDMD}. As
time progresses, we can update $\mA_k$ according to the formula given in
equation (\ref{eqn:windowA}).
%
%
However, evaluating~\eqref{eqn:windowA} involves computing a new pseudoinverse
and a matrix multiplication, which are costly operations.
We may compute this
update more efficiently using an approach similar to that in the previous
section, as we describe below.

\subsection{Algorithm for windowed DMD}
\label{sec:alg-window-dmd}
As in the approach presented in Section~\ref{sec:alg-online-dmd}, observe that equation (\ref{eqn:windowA}) can be written as
\begin{equation}
\label{eqn:windowdefAk}
\mA_k = \mY_k \mX_k^+ =  \mY_k \mX_k^T (\mX_k \mX_k^T)^{-1}=\mQ_k \mP_k,
\end{equation}
where
\begin{equation*}
\mQ_k = \mY_k \mX_k^T = \sum_{i=k-w+1}^{k} \vy_i \vx_i^T,
\end{equation*}
\begin{equation} \label{eqn:windowdefPk}
\mP_k = (\mX_k \mX_k^T)^{-1} = \bigg(\sum_{i=k-w+1}^{k} \vx_i \vx_i^T\bigg)^{-1}.
\end{equation}
where $\mQ_k$ and $\mP_k$ are $n\times n$ matrices. The condition that $\mX_k$
has rank~$n$ ensures that $\mX_k \mX_k^T$ is invertible, so $\mP_k$ is well defined.

At step $k+1$, we need to compute $\mA_{k+1} = \mQ_{k+1} \mP_{k+1}$. Clearly, $\mQ_{k+1}, \mP_{k+1}$ are related to $\mQ_{k}, \mP_{k}$. To show this, we write them down explicitly as
\begin{equation*}
\begin{split}
\mQ_{k+1} & = \mY_{k+1} \mX_{k+1}^T = \sum_{i=k-w+2}^{k+1} \vy_i \vx_i^T = \mQ_k -  \vy_{k-w+1} \vx_{k-w+1}^T + \vy_{k+1} \vx_{k+1}^T,
\end{split}
\end{equation*}
\begin{equation*}
\begin{split}
\mP_{k+1}^{-1}  & = \mX_{k+1} \mX_{k+1}^T = \sum_{i=k-w+2}^{k+1} \vx_i \vx_i^T =
\mP_k^{-1} -  \vx_{k-w+1} \vx_{k-w+1}^T + \vx_{k+1} \vx_{k+1}^T.
\end{split}
\end{equation*}
There is an intuitive interpretation to this relationship: $\mQ_{k}, \mP_{k}$
forgets the oldest snapshot and incorporates the newest snapshot, and gets updated
into $\mQ_{k+1}, \mP_{k+1}$.
As in Section~\ref{sec:alg-online-dmd}, where we used the Sherman-Morrison
formula~\eqref{eq:Sherman-Morrison} to update $\mP_k$, we may use a similar
approach to update $\mP_k$ in this case.

Letting
\begin{equation*}
\mU =
\begin{bmatrix}
\vx_{k-w+1} & \vx_{k+1}
\end{bmatrix},\qquad
\mV =
\begin{bmatrix}
\vy_{k-w+1} & \vy_{k+1}
\end{bmatrix},\qquad
\mC =
\begin{bmatrix}
-1 & 0 \\
0 & 1
\end{bmatrix},
\end{equation*}
we may write $\mQ_{k+1}, \mP_{k+1}$ as
\begin{equation*}
\mQ_{k+1} = \mQ_k + \mV \mC \mU^T,
\end{equation*}
\begin{equation*}
\mP_{k+1}^{-1} = \mP_k^{-1} + \mU \mC \mU^T,
\end{equation*}
therefore
\begin{equation}
\label{eqn:windowdefineAk1}
\mA_{k+1} = \mQ_{k+1} \mP_{k+1} = (\mQ_k + \mV \mC \mU^T)(\mP_k^{-1} + \mU \mC \mU^T)^{-1}.
\end{equation}
Now, the matrix inversion lemma (or Woodbury formula)
\cite{woodbury1950inverting, hager1989updating} states that
\begin{equation}
\label{eq:9}
(\mA + \mU \mC \mV)^{-1} = \mA^{-1} - \mA^{-1}\mU(\mC^{-1}+\mV \mA^{-1} \mU)^{-1}\mV \mA
\end{equation}
whenever $\mA$, $\mC$, and $\mA + \mU \mC \mV$ are invertible. Applying this
formula to our expression for $\mP_{k+1}$, we have
\begin{subequations}
\label{eq:7}
\begin{equation}
\label{eqn:windowupdatePk1}
\mP_{k+1}= \mP_k - \mP_k \mU \mGamma_{k+1} \mU^T \mP_k,
\end{equation}
where
\begin{equation}
\label{eqn:windowGamma}
\mGamma_{k+1} = (\mC^{-1}+\mU^T \mP_k \mU)^{-1}.
\end{equation}
\end{subequations}
Substituting back into~\eqref{eqn:windowdefineAk1}, we obtain
\begin{equation}
\label{eqn:windowexpandAk1}
\begin{split}
\mA_{k+1} & = (\mQ_k + \mV \mC \mU^T)(\mP_k - \mP_k \mU \mGamma_{k+1} \mU^T \mP_k) \\
& = \mQ_k\mP_k - \mQ_k \mP_k \mU \mGamma_{k+1} \mU^T \mP_k \\
& + \mV \mC \mU^T\mP_k - \mV \mC \mU^T \mP_k \mU \mGamma_{k+1} \mU^T \mP_k.
\end{split}
\end{equation}
The last two terms simplify, since
\begin{equation*}
\begin{split}
\mV \mC \mU^T\mP_k - \mV \mC \mU^T \mP_k \mU \mGamma_{k+1} \mU^T \mP_k & = \mV \mC (\mGamma_{k+1}^{-1} - \mU^T \mP_k \mU) \mGamma_{k+1} \mU^T \mP_k \\
& = \mV \mC \mC^{-1} \mGamma_{k+1} \mU^T \mP_k =  \mV  \mGamma_{k+1} \mU^T \mP_k,
\end{split}
\end{equation*}
where we have used~\eqref{eqn:windowGamma}. Substituting into~\eqref{eqn:windowexpandAk1}, we obtain
\begin{equation*}
\begin{split}
\mA_{k+1} & = \mQ_k\mP_k - \mQ_k \mP_k \mU \mGamma_{k+1} \mU^T \mP_k + \mV \mGamma_{k+1} \mU^T \mP_k\\
& = \mA_{k} - \mA_{k}\mU \mGamma_{k+1} \mU^T \mP_k + \mV \mGamma_{k+1} \mU^T \mP_k,
\end{split}
\end{equation*}
and hence
\begin{equation}
\label{eqn:windowupdateAk1}
\mA_{k+1} = \mA_{k} +(\mV -  \mA_{k}\mU) \mGamma_{k+1} \mU^T \mP_k.
\end{equation}
Notice the similarity between this expression with the updating
formula~\eqref{eqn:onlineupdateAk1} for online DMD. $\mGamma_{k+1}$ is the
matrix version of $\gamma_{k+1}$ in~\eqref{eqn:gammak1}. The matrix $(\mV -  \mA_{k}\mU)$ can also be considered as the prediction error based on current model $\mA_k$, and the correction to DMD matrix is proportional to this error term.

The updates in equations \eqref{eqn:windowupdateAk1}, \eqref{eq:7} require two
products of $n\times n$ and $n\times 2$ matrices (to compute $\mA_k \mU$ and
$\mP_k \mU$, since $\mP_k$ is symmetric), and two products of $n\times 2$
and $2\times n$ matrices, for a total of $8n^2$ multiplies.  This windowed DMD
approach is much more efficient than the standard DMD approach, solving~\eqref{eqn:windowdefAk} directly
($\mathcal{O}(wn^2)$ multiplies, with $w \geq n$).  Windowed DMD can be initialized
in the same manner as online DMD, discussed in Section~\ref{sec:alg-online-dmd}.

In order to implement windowed DMD, we need to store two $n \times n$ matrices
($\mA_k,\mP_k$), as well as the $w$ most recent snapshots.  Thus, the storage
required is more than in online DMD, or the weighted online DMD approach
discussed in Section~\ref{sec:weighted-online}, which also provides a mechanism
for ``forgetting'' older snapshots.

It is worth pointing out that the update formulas
\eqref{eq:7},\eqref{eqn:windowupdateAk1} give the exact solution $\mA_{k+1} =
\mY_{k+1} \mX_{k+1}^+$ from equation~(\ref{eqn:windowdefineAk1}), without approximation.

\paragraph{Larger window stride size}
We can in principle move more than one step for windowed DMD, i.e., forgetting
multiple snapshots and remembering multiple snapshots. If we would like to move
the sliding window for $s$ steps ($s<n/2$), then after similar derivations we
can show that the computational cost is $8sn^2$ multiplies, which is the same as
applying the rank-2 formulas $s$ times. Therefore, there is no obvious advantage
to incorporating multiple snapshots at one time.

\paragraph{Extensions}
Similar to online DMD, we can also incorporate an exponential weighting factor
into windowed DMD.  In particular, consider the cost function as
\begin{equation*}
\tilde J_k = \sum_{i=k-w+1}^k \rho^{k-i} \|\vy_i - \mA_k \vx_i\|^2,\qquad 0 < \rho \leq 1,
\end{equation*}
where $\rho$ is the weighting factor.  Then, proceeding as in
Section~\ref{sec:weighted-online}, we
obtain the update formulas
\begin{equation}
\label{eqn:windowupdateweightAk1}
\mA_{k+1} = \mA_{k} +(\mV -  \mA_{k}\mU) \tilde \mGamma_{k+1} \mU^T \hat \mP_k,
\end{equation}
\begin{subequations}
\begin{equation}
\hat\mP_{k+1} = \frac{1}{\rho} (\hat\mP_k - \hat\mP_k \mU \tilde \mGamma_{k+1} \mU^T \hat\mP_k ),
\end{equation}
where
\begin{equation}
\tilde \mGamma_{k+1} = (\tilde \mC^{-1}+\mU^T \hat\mP_k \mU)^{-1},\qquad
\tilde \mC =
\begin{bmatrix}
-\rho^w & 0 \\
0 & 1
\end{bmatrix}.
\end{equation}
\end{subequations}

As with the online DMD algorithm, the above windowed DMD algorithm applies
generally to extended DMD (EDMD) \cite{williams2015data} as well, if $\vx_k,
\vy_k$ are simply replaced by the observable vector of the states.  In addition,
the algorithm can be used for real-time online system identification, including
both linear and nonlinear system identification, as discussed in Section~\ref{sec:system-id}.

\paragraph{Summary} To summarize, the algorithm proceeds as follows:
\begin{enumerate}
\item Collect $w$ snapshot pairs $(\vx_j,\vy_j)$, $j=1,\ldots,w$, where $w \geq n$ is
  large enough so that $\rank \mX_k=n$ (where $\mX_k$ is given by~\eqref{eqn:windowXk}).
\item Compute $\mA_k$ and $\mP_k$ from~\eqref{eq:4}, where $\mX_k,\mY_k$ is given by~\eqref{eqn:windowXk}.
\item When a new snapshot pair $(\vx_{k+1},\vy_{k+1})$ becomes available, update
  $\mA_k$ and~$\mP_k$ according to~\eqref{eqn:windowupdateAk1} and~\eqref{eq:7}.

\end{enumerate}
Implementations of this algorithm in both Matlab and Python are publicly
available at~\cite{zhang2017odmd}.

\section{Online system identification}
\label{sec:system-id}
As previously mentioned, the online and windowed DMD algorithms discussed above
can be generalized to online system identification with control in a
straightforward manner.  For a review of system identification methods,
see~\cite{aastrom1971system}.



\subsection{Online linear system identification}
Dynamic Mode Decomposition can be used for system identification, as
shown in~\cite{proctor2016dynamic}.  Suppose we are interested in identifying a (discrete-time) linear system given by
\begin{equation}
\label{eq:8}
\vx_{k+1} = \mA \vx_k + \mB \vu_k,
\end{equation}
where $\vx_k \in \mathbb{R}^n, \vu_k \in \mathbb{R}^p$ are the states and control input respectively, $\mA \in \mathbb{R}^{n \times n}, \mB \in \mathbb{R}^{n \times p}$.

At time~$k$, assume that we have access to $\vx_1, \vx_2, \cdots, \vx_{k+1}$ and
$\vu_1, \vu_2, \cdots, \vu_{k}$.  Letting
\begin{equation*}
\tilde \mY_k=
\begin{bmatrix}
\vx_2 & \vx_3 &  \cdots & \vx_{k+1}
\end{bmatrix}, \qquad
\tilde \mX_k=
\begin{bmatrix}
\vx_1 & \vx_2 &  \cdots & \vx_{k}\\
\vu_1 & \vu_2 & \cdots & \vu_k
\end{bmatrix},\qquad
\tilde\mA=
\begin{bmatrix}
\mA & \mB
\end{bmatrix},
\end{equation*}
we may write~\eqref{eq:8} in the form
\begin{equation}
\tilde \mY_k = \tilde\mA \tilde\mX_k.
\end{equation}
The matrices $\mA,\mB$ may then be found by minimizing the cost function
\begin{equation}
J_k = \|\tilde \mY_k - \tilde\mA_k \tilde\mX_k\|_F^2.
\end{equation}
As before, the solution is given by
\begin{equation}
\tilde\mA_k = \tilde\mY_k \tilde\mX_k^+.
\end{equation}

At time $k+1$, we add a new column to $\tilde\mX_k$ and $\tilde\mY_k$, and we would like to
update $\tilde\mA_{k+1}$ using our previous knowledge of $\tilde\mA_k$.  Using
the same approach as in Section~\ref{sec:onlineDMD}, it is straightforward to
extend the online DMD and windowed DMD algorithms to this case. In particular,
the square matrix~$\mA_k$ from Section~\ref{sec:onlineDMD} is replaced by the
rectangular matrix $\tilde\mA_k$ defined above, and the vector $\vx_k$ in the
formulas in Section~\ref{sec:onlineDMD} is replaced by the column vector
\[
  \begin{bmatrix}
    \vx_k\\
    \vu_k
  \end{bmatrix}.
\]

\subsection{Online nonlinear system identification}
The efficient online/windowed DMD algorithms apply to nonlinear system
identification as well. In general, nonlinear system identification is a
challenging problem; see \cite{nelles2013nonlinear} for an overview. Some
interesting methods are to use linear-parameter-varying models
\cite{Verdult:2005,hemati2016parameter}, or to consider a large dictionary of potential
nonlinear functions, and exploit sparsity to select a small subset
\cite{brunton2016discovering}.

Suppose we are interested in identifying a
nonlinear system
\begin{equation*}
\vx_{k+1} = \bm{f}(\vx_k,\vu_k)
\end{equation*}
directly from data, where $\vx_k \in \mathbb{R}^n, \vu_k \in \mathbb{R}^p$ are
the state vector and control input respectively.
The specific form of nonlinearity is unknown, but in order to proceed, we have to make some assumptions about the nonlinear form. Assume that we have $q$ (nonlinear) observables $z_i(\vx,\vu),i = 1, 2, \cdots,q$, such that the underlying dynamics can be approximately described by
\begin{equation}
\label{eq:11}
\vx_{k+1} = \mA \bm{z}_k,
\end{equation}
where $\mA \in \mathbb{R}^{n \times q}$, and
\begin{equation*}
\bm{z}_k =
\begin{bmatrix}
z_1(\vx_k,\vu_k) & z_2(\vx_k,\vu_k) & \cdots & z_q(\vx_k,\vu_k)
\end{bmatrix}^T.
\end{equation*}
To illustrate how this representation works, we take $\vx \in \mathbb{R}, \vu \in \mathbb{R}$ for example, and assume the nonlinear dynamics is given by
\begin{equation*}
x_{k+1} = a_1 x_k + a_2 x_k^2 + a_3 u_k + a_4 u_k^2 + a_5 x_k u_k.
\end{equation*}
Then by setting $z_1(x,u) = x, z_2(x,u) = x^2, z_3(x,u) = u, z_4(x,u) = u^2, z_5(x,u) = xu$, we can write the dynamics in the form~\eqref{eq:11}, with
\begin{equation*}
\mA =
\begin{bmatrix}
a_1 & a_2 & a_3 & a_4 & a_5
\end{bmatrix}.
\end{equation*}
Note that in the above, the state $x_k$ still evolves nonlinearly (i.e.,
$x_{k+1}$ depends in a nonlinear way on $x_k$ and $u_k$), but we are able to
identify the coefficients~$a_k$ using linear regression (i.e., finding the
matrix~$\mA$ in~\eqref{eq:11}).
%

This approach is related to Carleman linearization~\cite{Bellman:1963}, although
in Carleman linearization, the goal is to find a true linear representation of
the dynamics in a higher-dimensional state space, and for most nonlinear
systems, it is not possible to obtain a finite-dimensional linear
representation.

In summary, by assuming a particular form of the nonlinearity, we can find the
coefficients of a nonlinear system using the same techniques as used in linear
system identification, writing the nonlinear system in the
form~\eqref{eq:11}.

\section{Application and results}
\label{sec:result}

In this section, we illustrate the methods on a number of examples, first
showing results for simple benchmark problems, and then using data from a wind
tunnel experiment.

\subsection{Benchmarks}
\label{sec:benchmarks}
We now present a study of the computational time of various DMD
algorithms. Two benchmark tasks are considered here.  In the first task, we wish
to know the DMD matrix only at the final time step, at which point we have
access to all of the data. In the second task, we wish to compute the DMD matrix
at each time, whenever a new snapshot is required.  The first task thus
represents the standard approach to computing the DMD matrix, while the second
task applies to situations where the system is time varying, and we
wish to update the DMD matrix in real time.

\paragraph{Asymptotic cost.}  First, we examine how the various algorithms scale with the state dimension~$n$
and the number of snapshots~$m$, for the two tasks described above.  In
particular, we are concerned with the over-constrained case in which $n<m$.  For
the standard algorithm, in which the DMD matrix is computed directly
using~\eqref{eqn:onlineA}, one must compute an $n\times m$ pseudoinverse and an
$n\times m$, $m\times n$ matrix multiplication.  For the first task, the
computational cost (measured by the number of multiplies) is thus
\begin{equation*}
T_{\text{standard}}=\mathcal{O}(nm \min(m,n)+mn^2) = \mathcal{O}(mn^2).
\end{equation*}
For the second task, in which we compute the DMD matrix at each time, we refer
to the standard algorithm as ``batch DMD'', since the snapshots are processed
all in one batch.  The method is initialized and applied after $m_0$ snapshots are
gathered (and in the examples below, we will take $m_0=n$), so the computational cost is
\begin{equation*}
T_{\text{batch}}=\mathcal{O}\bigg(\sum_{k=m_0}^{m}(nk \min(k,n)+kn^2)\bigg) = \mathcal{O}(m^2 n^2),
\end{equation*}
Next, we consider windowed DMD, for a window containing $w$ snapshots (with $n<w<m$).  In this
case, we refer
to the
standard algorithm, in which DMD matrix is computed directly
using~\eqref{eqn:windowA}, as ``mini-batch DMD''.  The
computational cost is given by
\begin{equation*}
T_{\text{mini-batch}}=\mathcal{O}\bigg(\sum_{k=w}^{m}(nw \min(n,w)+wn^2)\bigg) = \mathcal{O}(m w n^2),
\end{equation*}
For streaming DMD~\cite{hemati2014dynamic} for a fixed rank~$r$, the cost of one iteration is
$\mathcal{O}(r^2n)$, and for full-rank streaming DMD, the cost of one iteration
is $\mathcal{O}(n^2)$.  Thus, for either task, the overall cost after $m$
snapshots is
\begin{equation*}
T_{\text{streaming}}^{r=n}=\mathcal{O}(\sum_{k=1}^{m}n^2) = \mathcal{O}(m n^2),
\end{equation*}
and
\begin{equation*}
T_{\text{streaming}}^{r<n}=\mathcal{O}(\sum_{k=1}^{m}r^2n) = \mathcal{O}(m r^2 n).
\end{equation*}
(If, in streaming DMD, the compression step (step 3 in the algorithm
in~\cite{hemati2014dynamic}) is performed only every $r$ steps, then the cost is
reduced to $\mathcal{O}(mrn)$.)  Finally, for both online and windowed DMD algorithms, discussed in
Sections~\ref{sec:alg-online-dmd} and~\ref{sec:alg-window-dmd}, the cost per
timestep is $\mathcal{O}(n^2)$.  The algorithms are applied after
$w$ snapshots are gathered, so the overall cost of either algorithm is
\begin{equation*}
T_{\text{online}}=T_{\text{window}}=\mathcal{O}(\sum_{k=w+1}^{m}n^2) = \mathcal{O}(m n^2).
\end{equation*}

\paragraph{Results} We now compare the performance of the different algorithms
on actual examples, for the two tasks described above.  In particular, we
consider a system with state $\vx \in \mathbb{R}^n$, where $n$ varies between 2
and 1024.  The entries in the $n\times n$ matrix $\mA$ are chosen randomly,
according to a normal distribution (zero-mean with unit variance).  The
snapshots $\vx_1,\ldots,\vx_m$ are also chosen to be random vectors, whose
components are also chosen according to the standard normal distribution.  In
the tests below, we use a fixed number of snapshots $m=10^4$.  For mini-batch
DMD and windowed DMD, the window size is fixed at $w=2048$, and online DMD and
windowed DMD are both initialized after the first $w$ snapshot pairs.  For
streaming DMD with a fixed rank~$r$, we take $r=16$.  The simulations are
performed in MATLAB (R2016b) on a personal computer equipped with a 2.6 GHz
Intel Core i5 processor.

\begin{figure}
  \centering
  \begin{subfigure}[b]{0.5\textwidth}
    \begin{tikzonimage}[width=.9\linewidth]{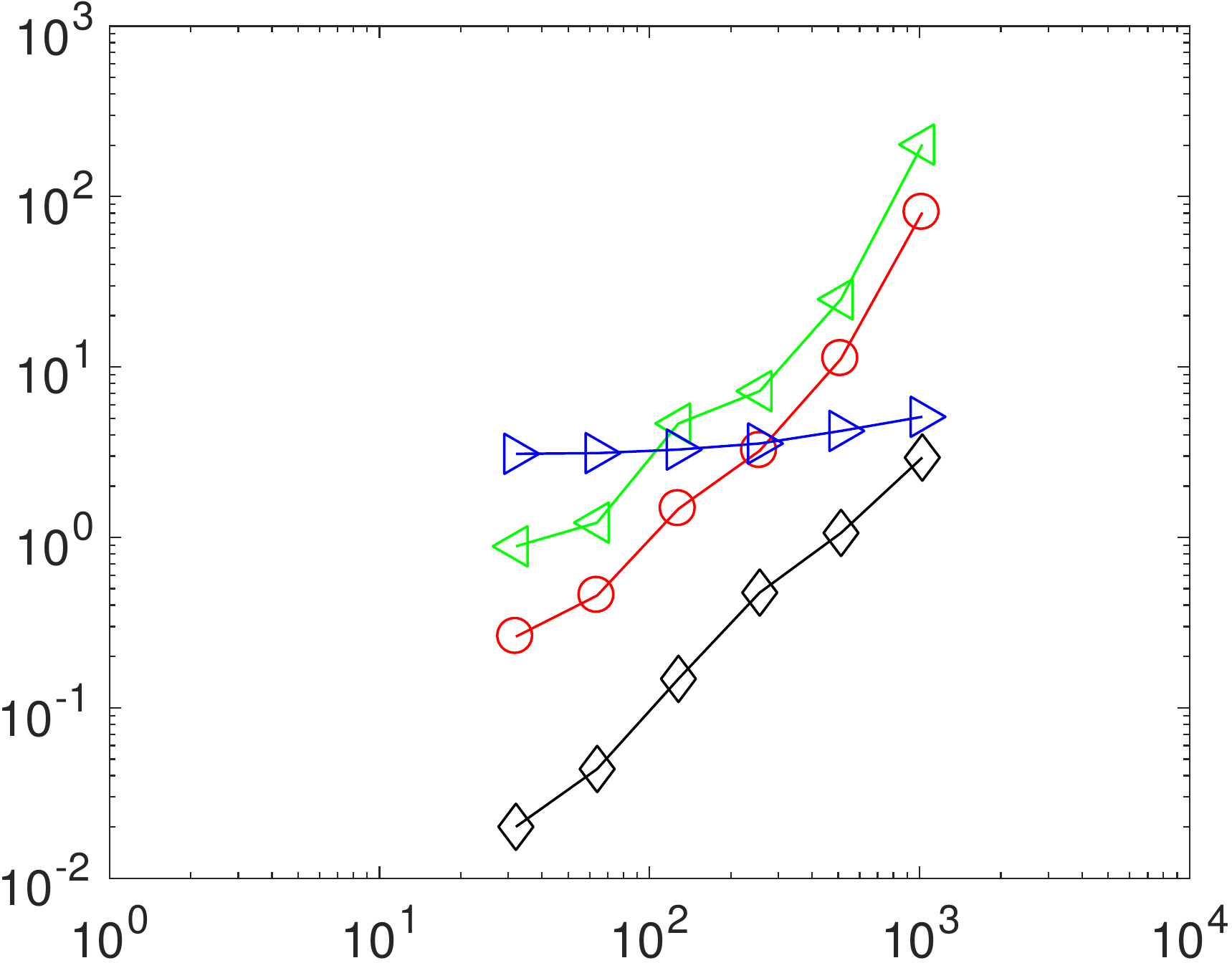}
      \small
      \node[below] at (0.5,0) {State dimension $n$};
      \node[anchor=south,rotate=90] at (0,0.5) {Time (sec)};
      \scriptsize
      \node[right] at (0.6,0.3) {Standard};
      \node[left, color=red] at (0.4,0.32) {Online};
      \node[left, color=green, align=right] at (0.4,0.43) {Streaming\\($\text{rank}=n$)};
      \node[left, color=blue, align=right] at (0.4,0.57) {Streaming\\($\text{rank}=r$)};
    \end{tikzonimage}
    \caption{Task: compute DMD matrix at final step.}
  \end{subfigure}%
  \begin{subfigure}[b]{0.5\textwidth}
    \begin{tikzonimage}[width=.9\linewidth]{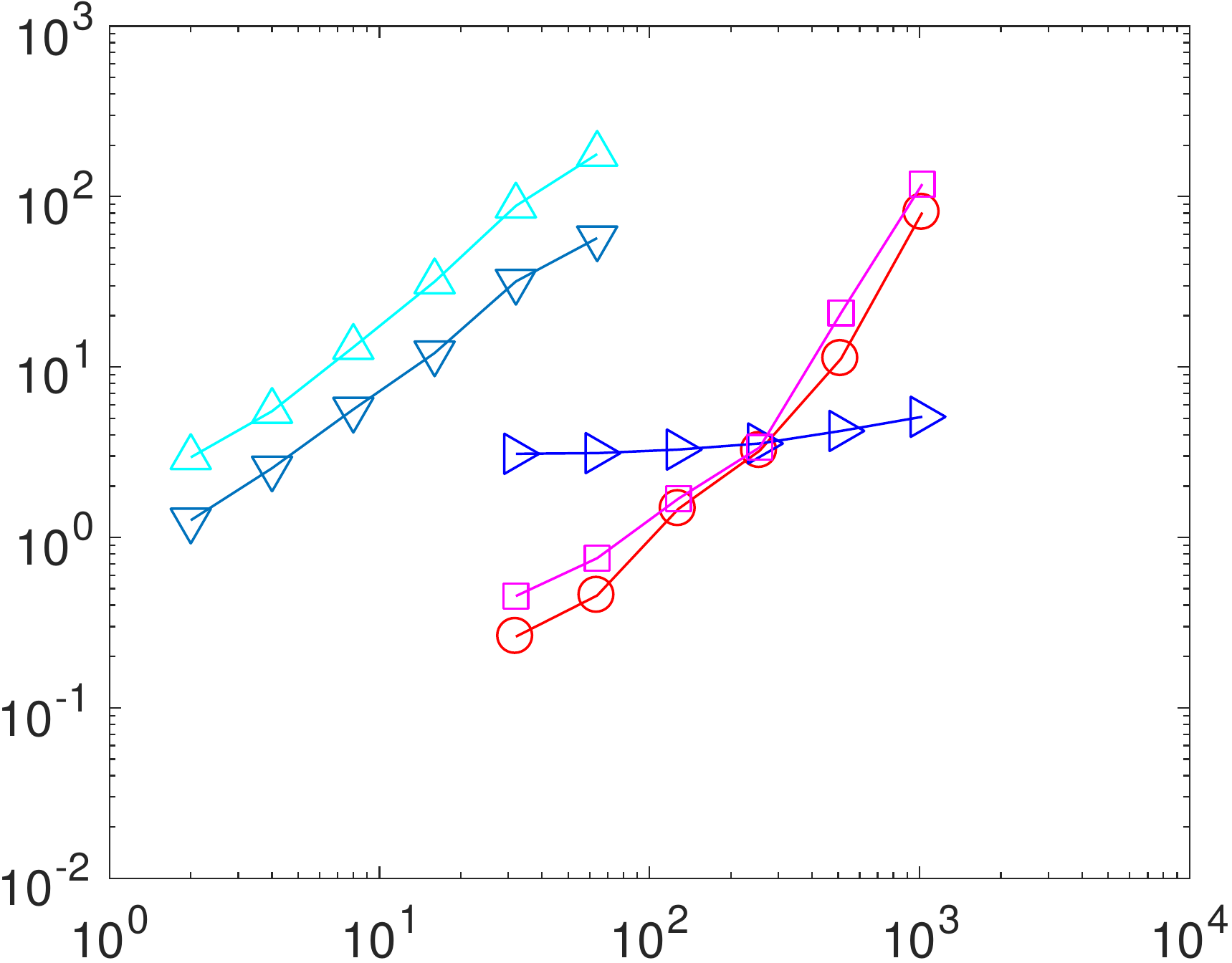}
      \small
      \node[below] at (0.5,0) {State dimension $n$};
      \node[anchor=south,rotate=90] at (0,0.5) {Time (sec)};
      \scriptsize
      \node[below, color=red] at (0.45,0.32) {Online};
      \node[left, color=magenta] at (0.4,0.39) {Windowed};
            \node[left, color=blue, align=right] at (0.95,0.47) {Streaming\\($\text{rank}=r$)};
      \node[right, color=blue!50] at (0.4,0.65) {Mini-batch};
      \node[left, color=cyan] at (0.38,0.78) {Batch};
    \end{tikzonimage}
    \caption{Task: compute DMD matrix at each step}
  \end{subfigure}%
  \caption{%
    Performance of the different DMD algorithms on the benchmark cases described
    in Section~\ref{sec:benchmarks}. For low-rank streaming, the dimension is
    limited to $r=16$.}
  \label{fig:time}
\end{figure}
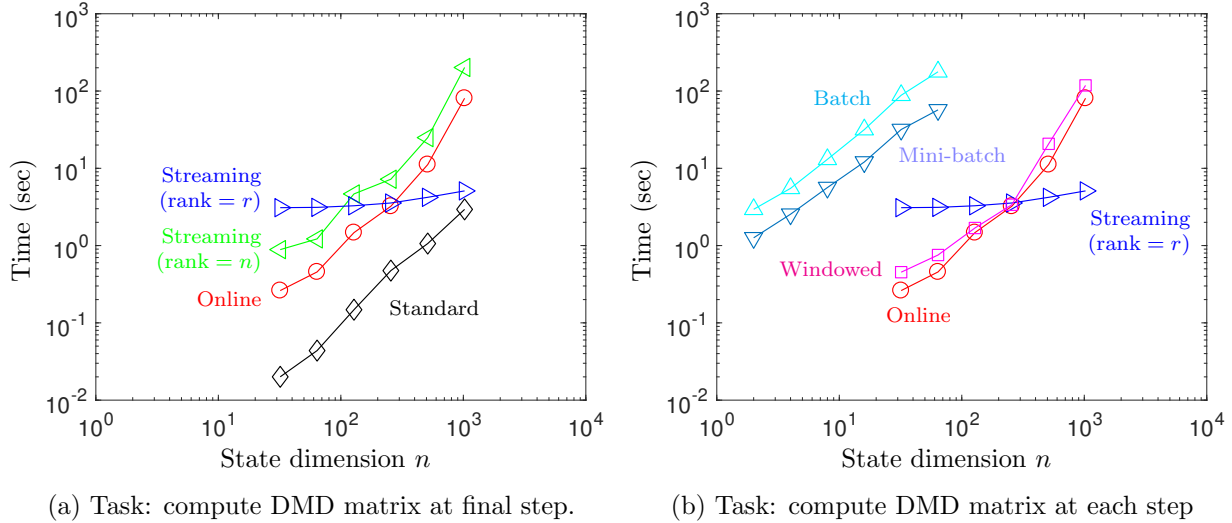

The results are shown in Figure~\ref{fig:time}.
For the first task (computing the DMD matrix only at the final step), the
standard DMD algorithm is the most efficient, for the problem sizes considered
here.  However, note that streaming DMD with a fixed rank~$r$ scales much better
with the state dimension~$n$, and would be the fastest approach for problems
with larger state dimension.

Our primary interest here is in the second task, shown in
Figure~\ref{fig:time}(b), in which the DMD matrix is updated at each step.  For
problems with $n<256$, online DMD is the fastest approach, and can be orders of
magnitude faster than the standard batch and mini-batch algorithms.  For
problems with larger state dimension, streaming DMD is the fastest algorithm, as
it scales linearly in the state dimension (while the other algorithms scale
quadratically).  However, note that streaming DMD does not compute the exact DMD
matrix: rather, it computes a projection onto a subspace of dimension $r$ (here
16).  By contrast, online DMD and windowed DMD both compute the full DMD matrix,
without approximation.

These results focus on the time required for these algorithms, but it is worth
pointing out the memory requirements as well.  Streaming DMD and online DMD do
not require storage of any past snapshots, while windowed DMD and mini-batch DMD
require storing the $w$ snapshots in the window, and batch DMD requires storage
of all past snapshots.

\subsection{Linear time-varying system}
\begin{figure}
  \centering
  \begin{subfigure}[b]{0.5\textwidth}
    \centering
    \begin{tikzonimage}[width=.8\linewidth]{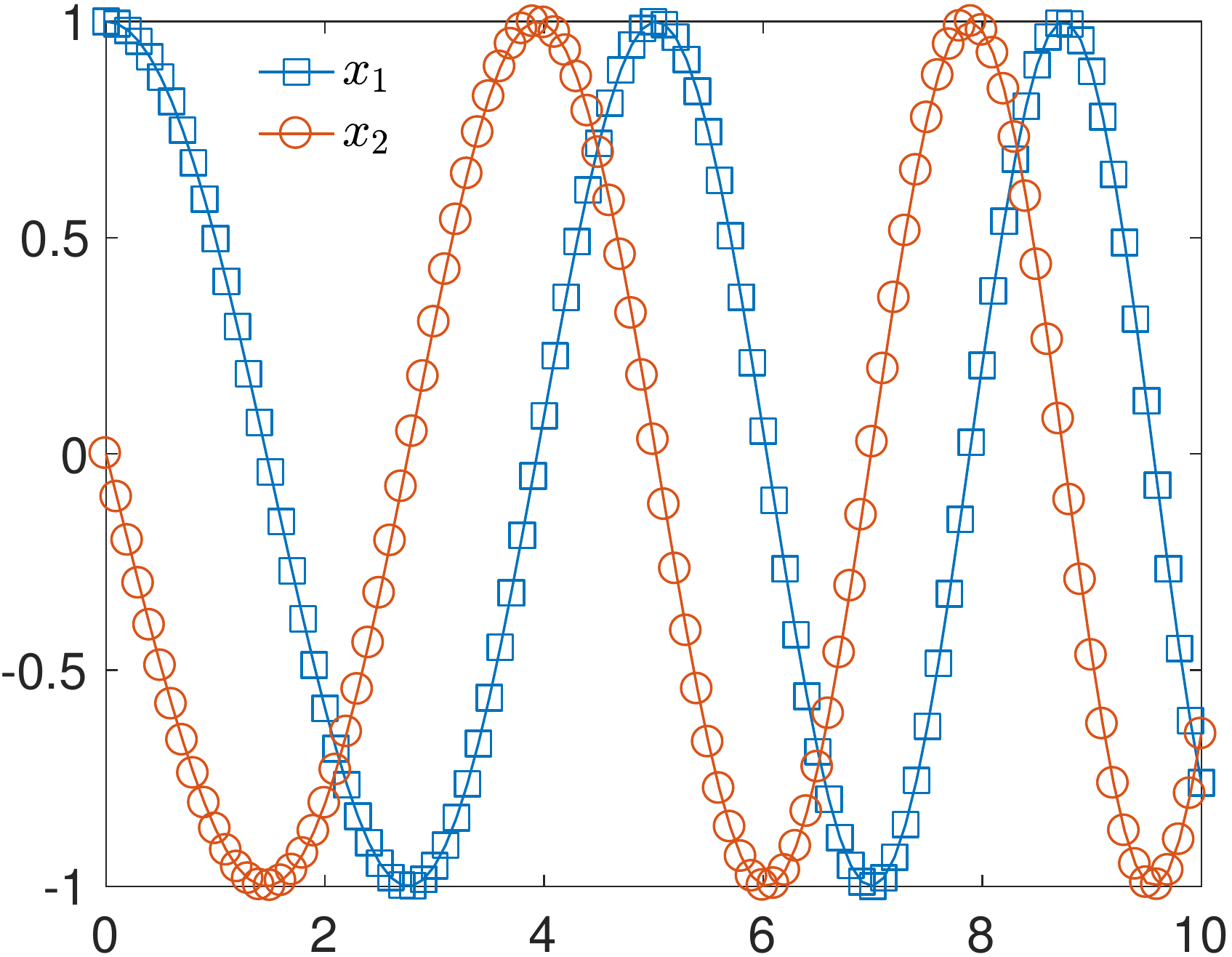}
      \small
      \node[below] at (0.5,0) {Time $t$};
      \node[anchor=south,rotate=90] at (0,0.5) {States $x_1,x_2$};
    \end{tikzonimage}
    \caption{State evolution}
  \end{subfigure}%
  \begin{subfigure}[b]{0.49\textwidth}
    \centering
    \begin{tikzonimage}[width=.8\linewidth]{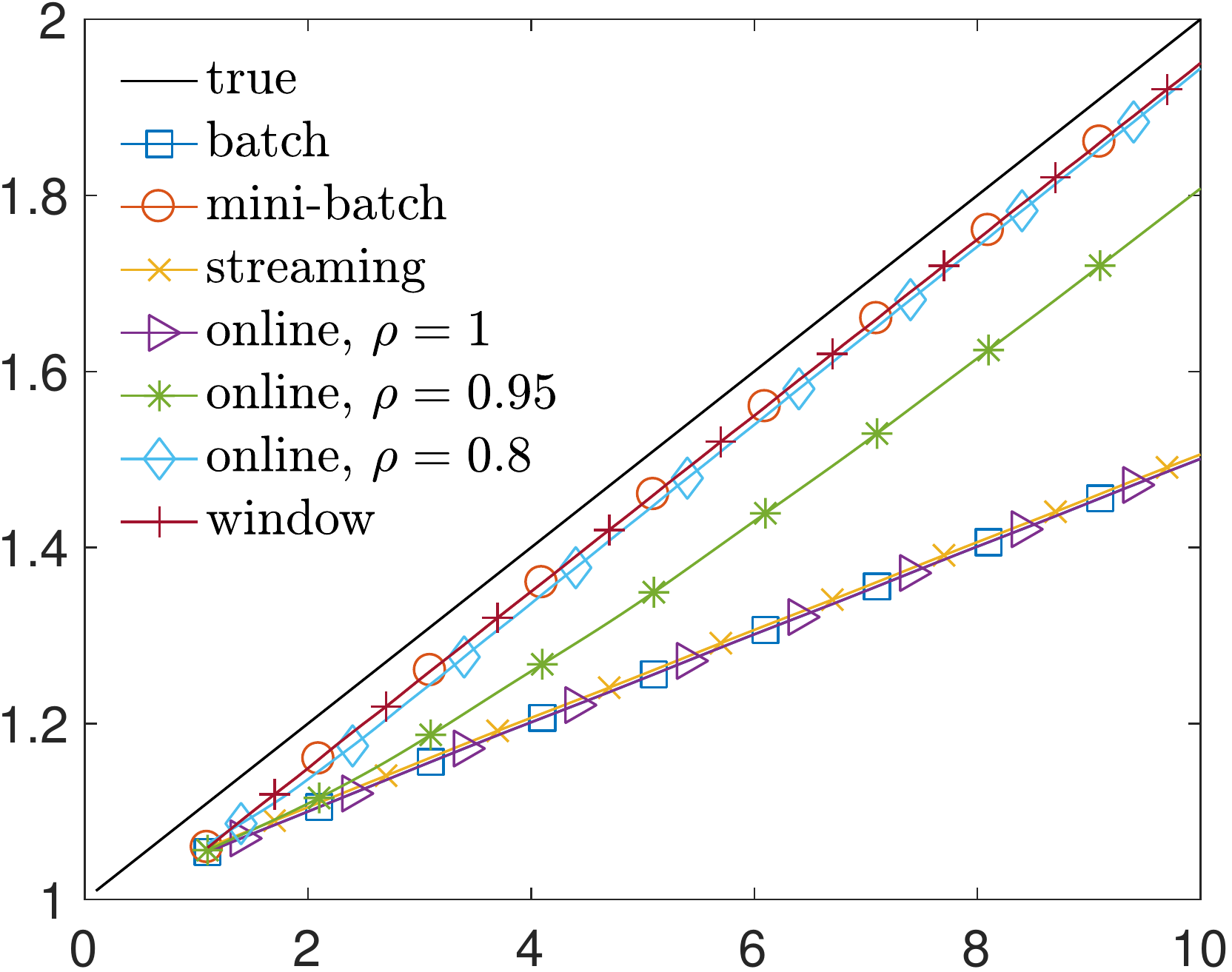}
      \small
      \node[below] at (0.5,0) {Time $t$};
      \node[anchor=south,rotate=90] at (0,0.5) {$\imag(\lambda)$};
    \end{tikzonimage}
    \caption{Frequency prediction}
  \end{subfigure}%
  \caption{%
    Solution of the linear time-varying system~\eqref{eq:ltv}, and frequencies
    predicted by various DMD algorithms.  For mini-batch and windowed DMD, the
    window size is $w=10$.  For online DMD, smaller values of the parameter
    $\rho$ result in faster tracking of the time-varying frequency.}
  \label{fig:2Dtimevarying}
\end{figure}

We now test the online DMD and windowed DMD algorithms on a simple linear
system that is slowly varying in time. In particular, consider the system
\begin{subequations}
\label{eq:ltv}
\begin{equation}
\dot{\vx}(t) = \mA(t) \vx(t),
\end{equation}
where $\vx(t) \in \mathbb{R}^2$, and the time-varying matrix $\mA(t)$ is given by
\begin{equation}
\mA(t) =
\begin{bmatrix}
0 & \omega(t) \\
-\omega(t) & 0
\end{bmatrix},
\end{equation}
\end{subequations}
where
\begin{equation*}
\omega(t)=1+\epsilon t.
\end{equation*}
We take $\epsilon=0.1$, so that the system is slowly varying in time. The
eigenvalues of $\mA(t)$ are $\pm i\omega(t)$, and it is straightforward to show
that $\|\vx(t)\|$ is constant in~$t$.  We simulate the system for $0<t<10$ from
initial condition $\vx(0)=(1,0)^T$, and the snapshots are taken with time step
$\Delta t=0.1$ as shown in Figure \ref{fig:2Dtimevarying}(a). It is evident from
the figure that the frequency is increasing with time.

Given the snapshots, we apply both brute-force batch DMD and mini-batch DMD as
benchmark, then we compare streaming DMD, online DMD and windowed DMD with these
two benchmark brute-force algorithms.
The finite time window size of mini-batch
DMD and windowed DMD is $w=10$. Batch DMD takes in account all the past
snapshots, while mini-batch DMD only takes the recent snapshots from a finite
time window. Streaming DMD, online DMD and windowed DMD are initialized using
the first $w=10$ snapshot pairs, and they start iteration from time $w+1$. Batch
DMD and mini-batch DMD also starts from time $w+1$. The results for streaming
DMD, online DMD ($\rho=1, \rho=0.95, \rho=0.8$), and windowed DMD are shown in
Figure \ref{fig:2Dtimevarying}(b). DMD finds
the discrete-time eigenvalues $\mu_{\text{DMD}}$ from data, and the
figure shows the continuous-time DMD eigenvalues $\lambda_\text{DMD}$, which are related to these by
\begin{equation}
\mu_\text{DMD} = e^{\lambda_\text{DMD}\Delta t},
\label{eqn:cevals}
\end{equation}
where $\Delta t$ is the time spacing between snapshot pairs.  We show the DMD
results starting from time $w+1$, and the true eigenvalues are also shown for comparison. 

Observe from Figure~\ref{fig:2Dtimevarying}(b) that the eigenvalues computed by
the standard algorithm (batch DMD) agree with those identified by streaming DMD
and online DMD with $\rho=1$, as expected.  Similarly, windowed DMD perfectly
overlaps with mini-batch DMD. When the weighting $\rho$ in online DMD is
smaller than~1, the identified frequencies shift slightly towards those
identified by windowed DMD. If we further decrease the weighting factor
($\rho=0.8$), online DMD aggressively forgets old data, and the identified
frequency adapts more quickly.  This example demonstrates that windowed DMD and
weighted online DMD are capable of capturing time-variations in dynamics, with
an appropriate choice of the weight~$\rho$.

\subsection{Pressure fluctuations in a separation bubble}
\label{sec:pressure-fluct}

We now demonstrate the algorithm on a more complicated example, using data
obtained from a wind tunnel experiment. In particular, we study the flow over a
flat plate with an adverse pressure gradient, and investigate the dynamics of
pressure fluctuations in the vicinity of a separation bubble.

\begin{figure}[!htb]
        \centering
     \includegraphics[width=.5\linewidth]{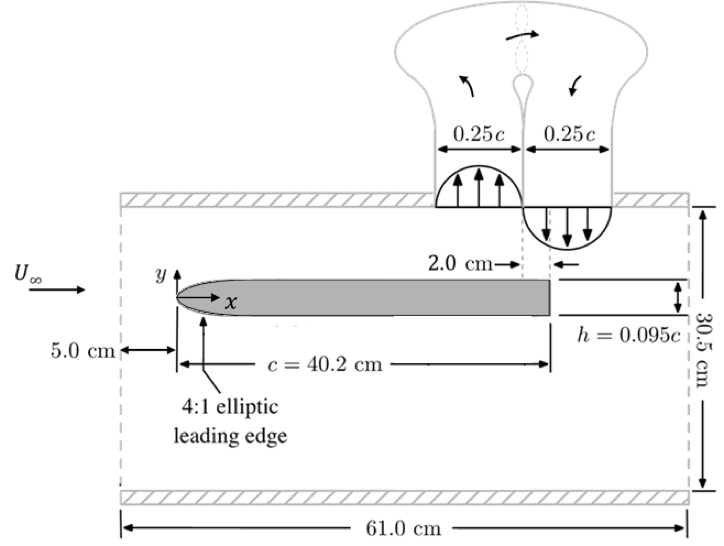}
     \caption{ Schematic of the flat plate model and flow separation system.}
\label{fig:config}
\end{figure}

The setup of the wind tunnel experiment is shown in Figure~\ref{fig:config}.  A
flat plate with a rounded leading edge is placed in the flow, and suction and
blowing are applied at the ceiling of the wind tunnel in order to apply a pressure
gradient at the surface of the plate, and cause the boundary layer to separate
and then re-attach.  The cross-section of the leading edge of the plate is a 4:1
ellipse, and the trailing edge of the model is square, which results in bluff-body shedding downstream.

These experiments were conducted in the Florida State Flow Control (FSFC)
open-return wind tunnel. The cross-sectional dimensions of the test section are
$30.5\,\text{cm}\times 30.5\,\text{cm}$, and the length is $61.0\,\text{cm}$.
The contraction ratio of the inlet is 9:1.  An aluminum honeycomb mesh and two
fine, anti-turbulence screens condition the flow at the inlet and provide a
freestream turbulence intensity of $u^\prime/U_{\infty} = 0.5\%$.  The
suction/blowing on the ceiling of the wind tunnel test section is provided by a
variable-speed fan mounted within a duct fixed to the ceiling of the wind
tunnel, which pulls flow from the ceiling and reintroduces it immediately
downstream.  The chord of the flat plate model is $c = 40.2$ cm, and the height
is $h = 0.095c$.  For these experiments, the freestream velocity is
$U_{\infty} = 3.9\,\text{m/s}$ and the Reynolds number is $\text{Re}_c =
U_\infty c/\nu = 10^5$.

Unsteady surface pressure fluctuations within the separated flow are monitored
by an array of 13 surface-mounted Panasonic WM-61A electret microphones located
within the separation region.  The microphones are placed at the centerline of
the plate, between $x/c=0.70$ and $0.94$, with a spacing of 0.02.  More details
regarding this microphone array can be found in \cite{DeemAIAA2017}.  These
13 microphone signals provide the data we use for online DMD.

Prior to applying online DMD, the microphone signals are conditioned to remove
external contaminating sources.  This process is described in
\cite{DeemAIAA2017}.  We collect pressure data for a total time of $T =
10\,\text{sec}$. The sampling rate of pressure snapshot is
$f_s=2048\,\text{Hz}$, so the time spacing between pressure snapshot is
$\Delta t = 1/f_s$, and the total number of snapshots is~$m=20480$. The state dimension is $n=13$, because there are 13 pressure sensors.

\begin{figure}[!htb]
  \centering
  \begin{tikzonimage}[width=.5\linewidth]{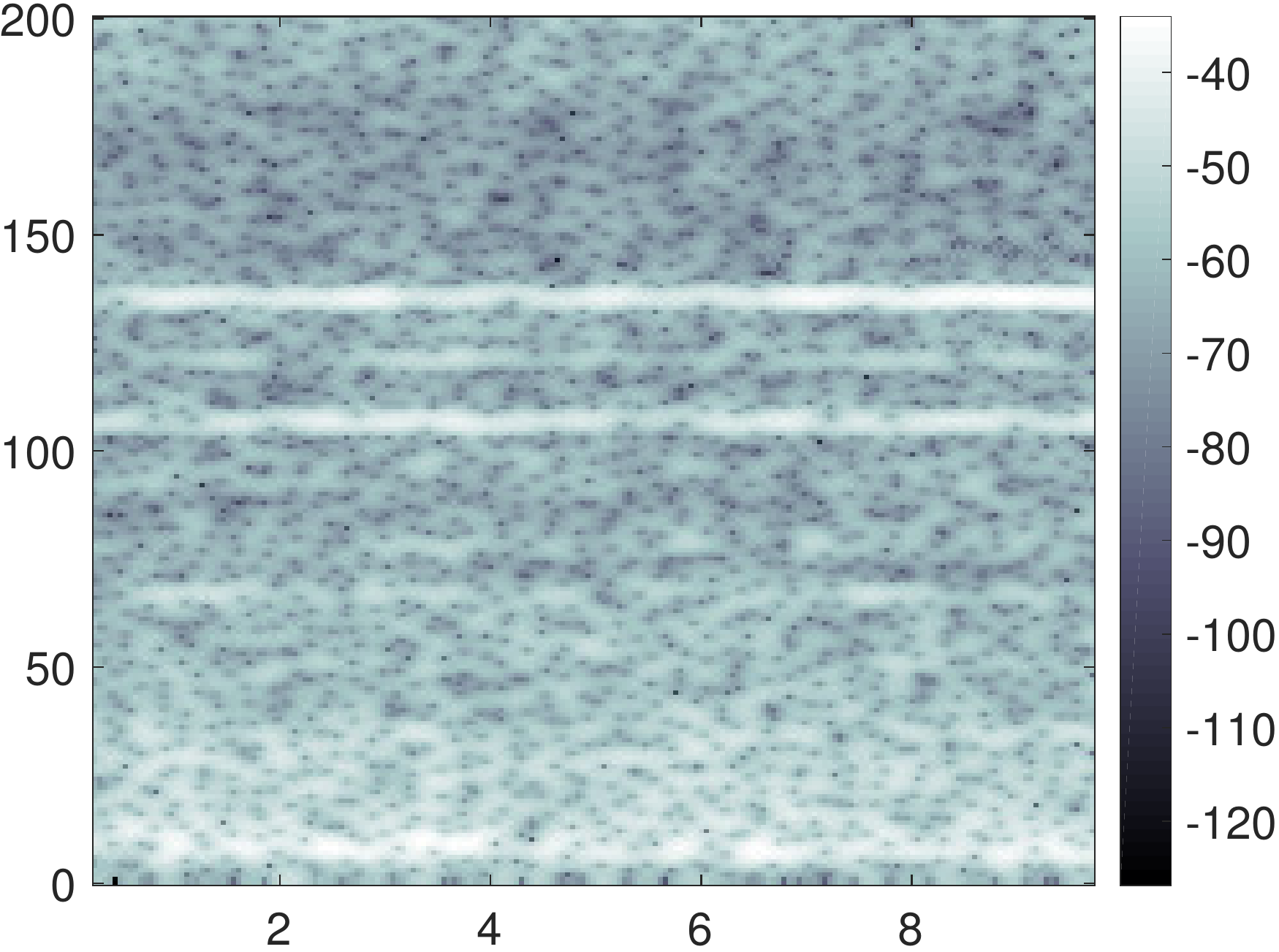}
    \small
    \node[below] at (0.5,0) {Time (sec)};
    \node[anchor=south,rotate=90] at (0,0.5) {Frequency (Hz)};
    \node[right] at (1,0.5) {PSD (dB/Hz)};
  \end{tikzonimage}
  \caption{%
    Power spectral density (PSD) of the pressure measurement at the first
    (upstream) pressure
    sensor.  Note the dominant frequencies at 105\,Hz and 135\,Hz.}
  \label{fig:spectrogram}
\end{figure}

We first present a spectral analysis of the pressure data using short-time
discrete Fourier transform (DFT). Figure~\ref{fig:spectrogram} shows the results
for the first (upstream) pressure sensor; other pressure sensors have similar
results. A window size of $w=1000$ is used, with overlap of $900$ samples between adjointing sections.

It is observed that two dominant frequencies (at about 105\,Hz and 135\,Hz) are
present over the whole time interval, while fluctuations at other frequencies
are slowly varying with time. We may use DMD to gain a comprehensive
understanding of the frequency variations in all the pressure sensors, and how
these might be related to one another.

\begin{figure}[!htb]
  \centering
  \begin{tikzonimage}[width=.9\linewidth]{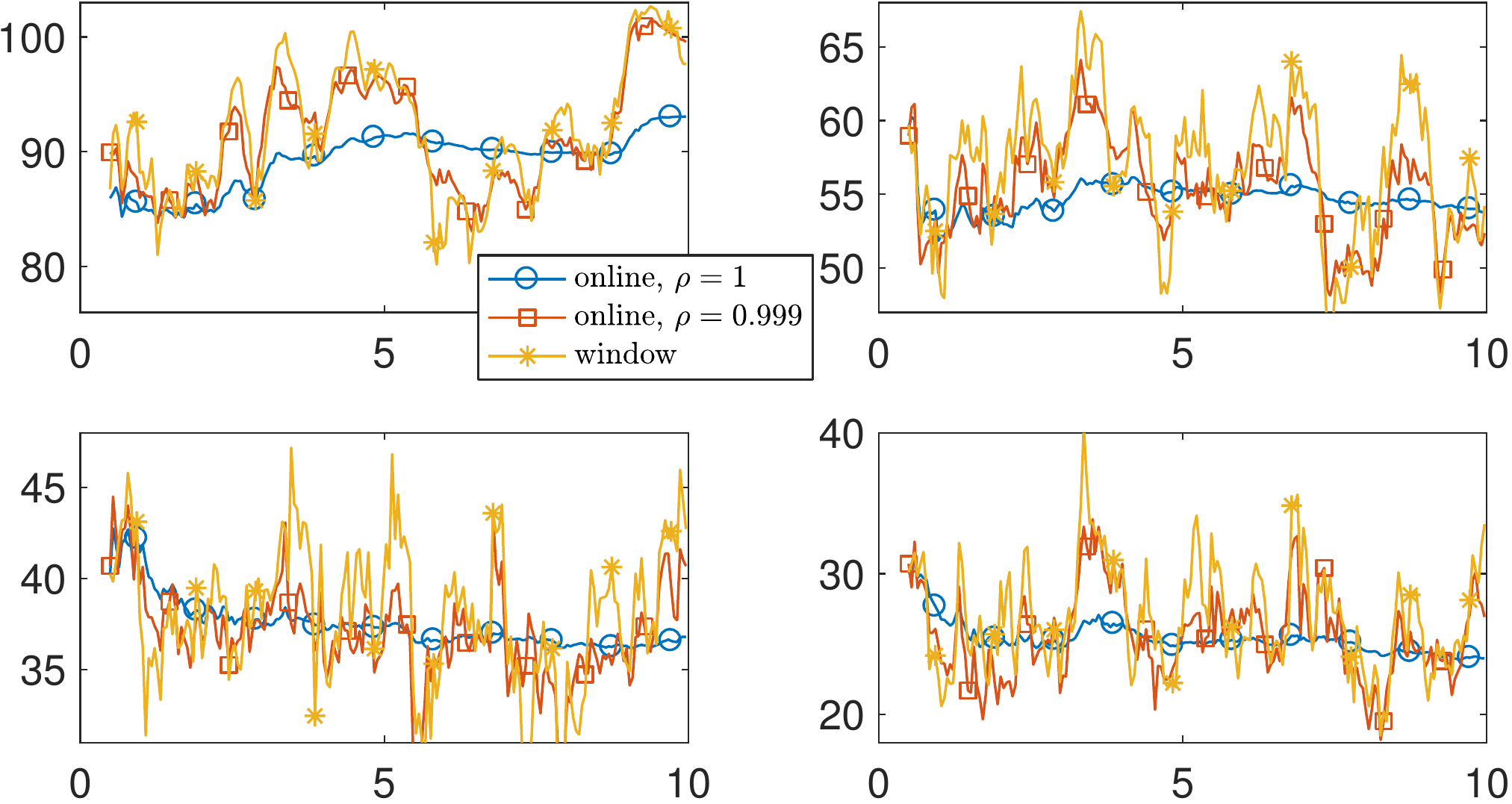}
    \small
    \node[below] at (0.25,0) {Time (sec)};
    \node[below] at (0.78,0) {Time (sec)};
    \node[anchor=south,rotate=90] at (0,0.25) {Frequency (Hz)};
    \node[anchor=south,rotate=90] at (0,0.8) {Frequency (Hz)};
  \end{tikzonimage}
  \caption{Four dominant DMD frequencies identified by different DMD algorithms
    from 13 pressure signals, as described in Section~\ref{sec:pressure-fluct}.}
  \label{fig:pressureDMD}
\end{figure}
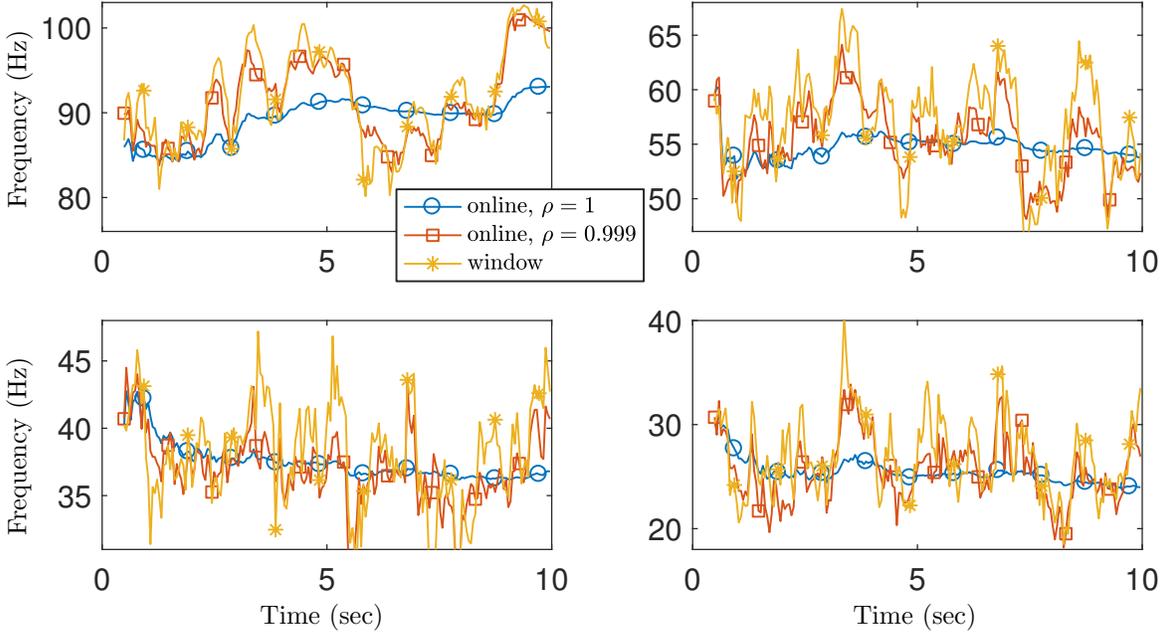

Next, we apply online DMD and windowed DMD to the
pressure dataset obtained from the experiment. Observe that the number of
snapshots $m=20480$ is much larger than the state dimension is $n=13$, so the
over-constrained assumption is satisfied.
The dynamics of the pressure fluctuations can be characterized by the DMD
frequencies, which may be slowly varying in time.
The DMD frequency is defined as
$$f_{\text{DMD}}=\frac{\imag(\lambda_{\text{DMD}})}{2\pi},$$
where $\imag(\lambda_{\text{DMD}})$ is the imaginary part of the continuous-time
DMD eigenvalues computed from equation~\eqref{eqn:cevals}. (The discrete-time
eigenvalues $\mu_\text{DMD}$ are eigenvalues of the $13\times 13$ matrix
matrix~$\mA_k$.) The four dominant frequencies computed by the various DMD algorithms are shown in
Figure~\ref{fig:pressureDMD}.
There are 13 DMD eigenvalues in total, and one of them is $f_0=0$ corresponding
to the mean flow. The remaining DMD eigenvalues consist of six complex-conjugate
pairs, corresponding to six non-zero DMD frequencies. For visualization, we show
only the four most dominant DMD frequencies, starting from time step $w+1$.  For
windowed DMD, we use a window size of $w=1000$, and for weighted online DMD, we
use $\rho=0.999$.

Recall that with $\rho=1$, online DMD coincides with the standard DMD
algorithm.  From Figure~\ref{fig:pressureDMD}, we see that for this case, the
frequencies remain more or less constant in time.  With $\rho=0.999$, online DMD
behaves more like windowed DMD: in particular, the method is better at tracking
variations in the frequency.
For online DMD with $\rho=0.999$, note that snapshot 1000 (the last included in windowed DMD) is
given a weight $0.999^{1000}\approx 0.37$.  The weighting factor in online DMD
acts like a soft cut-off for the old snapshots, compared with the hard cut-off
imposed by windowed DMD.
While the frequency variations shown in Figure~\ref{fig:pressureDMD}
appear to be rapid or noisy, note that the time interval shown in the figure is quite
long (about 300 periods for the lowest frequency of around 30\,Hz), so it is
reasonable to consider these frequencies as slowly varying in time.

\section{Conclusion and outlook}
In this work, we have developed efficient methods for computing online DMD and
windowed DMD. The proposed algorithms are especially useful in applications for
which the number of snapshots is very large compared to the state dimension, or
when the dynamics are slowly varying in time.  A weighting factor can be
included easily in the online DMD algorithm, which is used to weight recent
snapshots more heavily than older snapshots.  This approach corresponds to using
a soft cutoff for older snapshots, while windowed DMD uses a hard cutoff, from
a finite time window. The proposed algorithms can be readily extended to online
system identification, even for time-varying systems.

\begin{table}[!htb]
\centering
\renewcommand{\arraystretch}{1.2}
\begin{tabular}{ |c|c|c|c|c|c|c| }
 \hline
 Aspect & Standard & Batch & Mini-batch & Streaming & Online & Windowed\\
  \hline \hline
 Computational time & $\mathcal{O}(mn^2)$ & $\mathcal{O}(kn^2)$ & $\mathcal{O}(wn^2)$ & $\mathcal{O}(r^2n)$ & $4n^2$ & $8n^2$ \\
 Memory & $mn$ & $kn$ & $wn$ & $\mathcal{O}(rn)$ & $2n^2$ & $wn+2n^2$  \\
 Store past snapshots & Yes & Yes & Yes & No & No & Yes  \\
Track time variations & No & No & Yes & Yes & Yes & Yes  \\
Real-time DMD matrix & No & Yes & Yes & Yes & Yes & Yes  \\
Exact DMD matrix & Yes & Yes & Yes & No & Yes & Yes  \\
 \hline
\end{tabular}
\caption{Characteristics of the various DMD algorithms considered. Relevant parameters are state dimension $n$, total number of snapshot pairs $m \gg n$, window size $w$ such that $n<w \ll m$, low rank $r<n$, and discrete time $k>n$. Computational time denotes the required floating-point multplies for one iteration (computing the DMD matrix).}
\label{table:DMD}
\end{table}

The efficiency is compared against the standard DMD algorithm, both for
situations in which one computes the DMD matrix only at the final time, and for
situations in which one computes the DMD matrix in an ``online'' manner,
updating it as new snapshots become available.  The latter case is applicable,
for instance, when one expects the dynamics to be time varying.  For the former
case, the standard DMD algorithm is the most efficient, while for the latter
case, the new online and windowed DMD algorithms are the most efficient, and can
be orders of magnitude more efficient than the standard DMD algorithm.  Table
\ref{table:DMD} provides a brief comparison of the main characteristics and
features of standard DMD, batch DMD, mini-batch DMD, (low rank) streaming DMD,
online DMD, and windowed DMD.

The algorithms are further demonstrated on a number of examples, including a
linear time-varying system, and data obtained from a wind tunnel experiment.  As
expected, weighted online DMD and windowed DMD are effective at capturing
time-varying dynamics.

A straightforward and relevant direction for future work is more detailed study
of the application of proposed online/windowed DMD algorithms to system
identification. In cases where there are variations in dynamics, or where we
desire real-time control, it is crucial to build accurate and adaptive reduced
order models for effective control, and the methods proposed here could be
useful in these cases.

\section*{Acknowledgment}
We gratefully acknowledge funding from the Air Force Office of Scientific Research (AFOSR) grant FA9550-14-1-0289, monitored by Dr.~Doug Smith, and by DARPA award HR0011-16-C-0116.

\bibliography{ref}
\end{document}